\documentclass{article}%
\usepackage{amsmath}
\usepackage{amsfonts}
\usepackage{amssymb}
\usepackage{graphicx}
\usepackage{setspace}
\usepackage{hyperref}
\usepackage{mathtools}
\DeclarePairedDelimiter\floor{\lfloor}{\rfloor}
\usepackage{latexsym}
\usepackage{psfrag}
\usepackage{subfigure}

\hypersetup{colorlinks=true, linkcolor=blue, citecolor=red, urlcolor=red}
\newtheorem{theorem}{Theorem}

\newtheorem{lemma}[theorem]{Lemma}

\newenvironment{proof}[1][Proof]{\noindent\textbf{#1.} }{\ \rule{0.5em}{0.5em}}

\begin{document}
\doublespacing
\title{\textbf{Kullback-Leibler Divergence for Bayesian Nonparametric Model Checking}}   
\author{Luai Al-Labadi\thanks{{\em Address for correspondence:} Luai Al-Labadi, Department  of Mathematics, University of Sharjah, P. O. Box 27272, Sharjah, UAE. E-mail:  lallabadi@sharjah.ac.ae}, Viskakh Patel \thanks{ Department of Mathematical \& Computational Sciences, University of Toronto Mississauga, Ontario L5L 1C6, Canada. E-mail: vishakh.patel@mail.utoronto.ca}, Kasra Vakiloroayaei \thanks{ Department of Mathematical \& Computational Sciences, University of Toronto Mississauga, Ontario L5L 1C6, Canada. E-mail: k.vakiloroayaei@mail.utoronto.ca},  Clement Wan \thanks{ Department of Mathematical \& Computational Sciences, University of Toronto Mississauga, Ontario L5L 1C6, Canada. E-mail: cmclement.wan@mail.utoronto.ca}}
%


\date{\vspace{-5ex}}
\maketitle
\pagestyle {myheadings} \markboth {} {KL Divergence for Bayesian Nonparametric Model Checking}
\begin{abstract}
Bayesian nonparametric statistics is an area of considerable  research interest. While recently  there has been an extensive concentration in developing Bayesian nonparametric procedures for model checking, the use of the Dirichlet process, in its simplest form, along with the Kullback-Leibler divergence is  still an open problem.   This is mainly attributed to the discreteness property of the Dirichlet process and that the Kullback-Leibler divergence between any discrete distribution and any continuous distribution is infinity. The approach proposed in this paper, which is  based on  incorporating the Dirichlet process, the Kullback-Leibler divergence and the relative belief ratio, is considered the first concrete solution to this issue. Applying the approach is simple and does not require obtaining a closed form of the relative belief ratio.  A Monte Carlo study and real data examples show that the developed approach exhibits excellent performance.

\par
 \vspace{9pt} \noindent\textsc{Keywords:}  Bayesian Non-parametric, Dirichlet process,  Kullback-Leibler divergence, Model checking, Relative belief ratio.

 \vspace{9pt}

\end{abstract}

\section{Introduction}
Let $x=(x_{1},\ldots,x_{n})$ be a sample from a distribution $P$. The goal is to assess the hypothesis $\mathcal{H}_{0}:P\in\left\{  F_{\theta}:\theta
\in\Theta\right\}  $, where  $\left\{  F_{\theta}:\theta\in\Theta\right\}  $ denotes the collection of continuous cumulative distribution functions (cdf's). This problem is known as \emph{model checking} and it is quiet important in statistics. For instance, Jordan (2011) placed model checking and hypothesis testing as number one in a list of top-five open problems in Bayesian statistics.

Several Bayesian nonparametric procedures have been developed for model cehing.
A main approach considers  embedding the proposed model as a null hypothesis
in a larger family of distributions. Then priors are placed on the null and
the alternative and a Bayes factor is computed. Using a Dirichlet process for the prior on the alternative can be found by Carota and Parmigiani (1996), and Florens, Richard, and Rolin(1996). Verdinelli and Wasserman (1998),
Berger and Guglielmi (2001) and McVinish, Rousseau, and Mengersen (2009)
considered other types of  priors on the
alternative. Another  important approach utilized for model testing is to place a prior on the true distribution that is generating the data and then measuring the distance between the posterior distribution and the hypothesized one. Swartz (1999) and Al-Labadi and Zarepour (2013, 2014) used the Dirichlet process as a prior and then considered the Kolmogorov distance in order to derive a goodness-of-fit test for continuous models. To test for  discrete models, Viele (2007) used the Dirichlet process and the Kullback-Leibler (KL) divergence. For continuous model, Viele commented that his method ``cannot be used for continuous data directly because the Dirichlet Process is discrete with probability 1. The KL information between any discrete
distribution and any continous distribution is infinity, and thus we must find a
nonparametric method that produces continuous distributions. We employ a Dirichlet
Process Mixture (DPM)." In fact working with the Dirichlet Process Mixture adds some complexity to the approach and  makes it hard to implement   by many users. Hsieh (2011) used the P\'olya tree as a prior and the Kullback-Leibler divergence to test  for continuous distributions. To judge whether a resulting divergence measure is large or
small, he used  normal approximations based on  running a regression
of the means and standard deviations. Al-Labadi and Evans (2018) established a new approach for model checking by utilizing the Dirichlet process and relative belief ratios. Then to measure the change from a priori to a posteriori they used Cram\'er-von Mises distance. See also   Al-Labadi (2018) and  Al-Labadi, Zeynep and Evans (2017, 2018) and Evans and Tomal (2018) for examples of using relative belief ratios in different hypothesis testing problems.

Although the KL divergence sits atop most distance/divergence measures (Viele, 2007), it follows clearly from the previous discussion that its use alongside the Dirichlet process is very limited.  This is mainly due to the discreteness property of the  Dirichlet process.  A complete solution  to this obstacle is offered throughout this paper, where the Dirichlet process  is considered as a prior on
$P$ (the true/sampling distribution). Then the  concentrations of the distribution of the KL divergence between the prior and the model of interest is compared to that between the posterior and the model.  If the posterior is more concentrated about the model than the prior, then this is evidence in favor of the model and if the posterior is
less concentrated, then this is evidence against the model. See Figure 1 below, which represents  a plot of the prior and posterior densities of the KL divergence  when $\mathcal{H}_{0}$ is true and indeed the posterior is much more concentrated about 0 than the prior.  The comparison is
made via a relative belief ratio (Evans, 2015), which measures the evidence in the observed data
for or against the model, and a measure of the strength of this evidence is
also provided; so the methodology is based on a direct measure of statistical
evidence. Implementing the approach is direct and does not require
obtaining a closed form of the relative belief ratio. In addition, the  methodology does not require the use of a prior on $\theta$ and so is truly a check on the model itself avoiding any issues with the prior on $\theta.$

\begin{figure}[h]
\centering
\includegraphics[width=8cm]{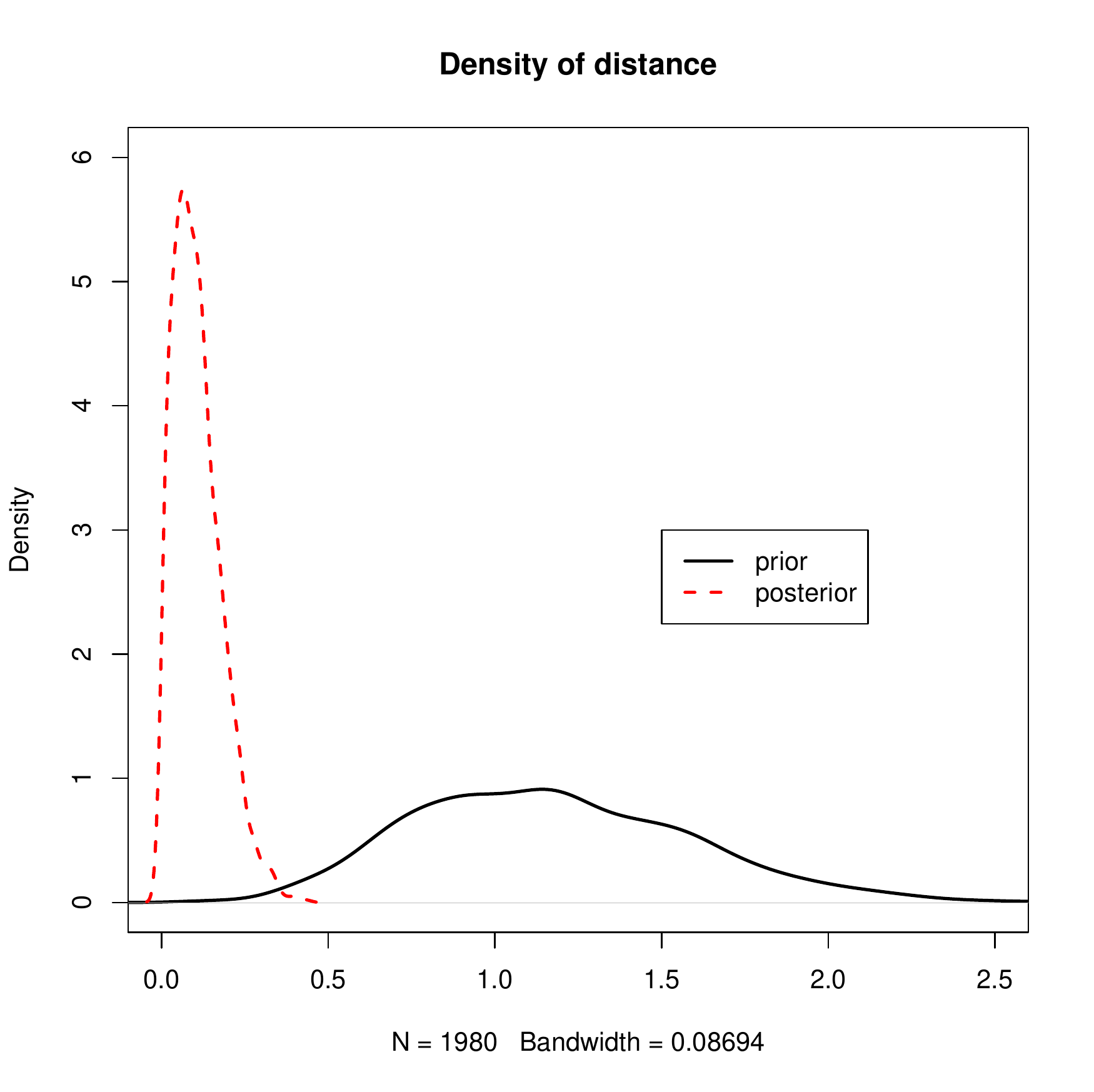}
\caption{Plot of prior density versus posterior density of distance when the model is correct (Example 1: $P_{true}=N(0,1)$). Clearly, the posterior distribution of the KL distance is more concentrated about 0 than that of the prior distance.}%
\end{figure}

This paper is organized as follows. In Section 2 and Section 3, the relative belief ratio and the Dirichlet process are briefly reviewed, respectively.  In Section 4, the Kullback-Leibler divergence between probability measures is discussed. Section 5 discusses the proposed approach for model checking, where it is argued that a particular choice of
the Kullback-Leibler divergence and the Dirichlet process should be employed. In Section 6, a computational algorithm for the implementation of the approach is outlined. Section 7 presents a
number of examples where the behavior of the methodology is examined in some
detail. Section 8 ends with a brief summary of the  results.

\section{Relative Belief Ratios}

Consider $\{f_{\theta}:\theta\in\Theta\}$ to be a collection of densities on a
sample space $\mathcal{X}$ and let $\pi$ be a prior on $\Theta.$ Given the data $x$, the posterior distribution of $\theta$ is  $\pi(\theta\,|\,x)=\pi(\theta)f_{\theta}(x)/\int
_{\Theta}\pi(\theta)f_{\theta}(x)\,d\theta$. Let $\psi=\Psi(\theta)$  be the parameter of interest. Then the
prior and posterior densities of $\psi$ are denoted by $\pi_{\Psi}$ and $\pi_{\Psi}%
(\cdot\,|\,x),$ respectively. The relative belief ratio (Evans, 2015) for a value $\psi$ is
then defined as  $RB_{\Psi}(\psi\,|\,x)=\lim_{\delta\rightarrow0}\Pi_{\Psi
}(N_{\delta}(\psi\,)|\,x)/\Pi_{\Psi}(N_{\delta}(\psi\,))$, where $N_{\delta
}(\psi\,)$ is a sequence of neighbourhoods of $\psi$ converging   nicely (see, for example, Rudin (1974)) to
$\psi$ as $\delta\rightarrow0.$ More commonly,
\begin{equation}
RB_{\Psi}(\psi\,|\,x)=\pi_{\Psi}(\psi\,|\,x)/\pi_{\Psi}(\psi), \label{relbel}%
\end{equation}
the ratio of the posterior density to the prior density at $\psi.$ That is,
$RB_{\Psi}(\psi\,|\,x)$ is measuring how beliefs have changed that
$\psi$ is the true value from \textit{a priori} to \textit{a posteriori}. Note that, a relative
belief ratio is similar to a Bayes factor, as both are measures of evidence,
but the latter measures this via the change in an odds ratio.  A discussion about the relationship between relative belief ratios and Bayes factors is detailed in Baskurt and Evans (2013). In particular, when a Bayes factor is defined via a limit in the continuous
case, the limiting value is the corresponding relative belief ratio.

By a basic principle of evidence, $RB_{\Psi}(\psi\,|\,x)>1$ implies that the probability of $\psi$ being correct increases after observing the data, and so there is evidence in favour of $\psi$. Else if $RB_{\Psi}(\psi\,|\,x)<1$ then the data claims evidence of the $\psi$ being incorrect and thus evidence of against $\psi$. Also if the $RB_{\Psi}(\psi\,|\,x)=1$, then there is no evidence either way.

Therefore, the $RB_{\Psi}(\psi_{0}\,|\,x)$ measures the evidence of the hypothesis ${H}_{0}=\{\theta:\Psi(\theta)=\psi_{0}\}.$ It is critical to rectify the degree of strength and weakness of this value. One of nicer calibration of the $RB_{\Psi}(\psi_{0}\,|\,x)$ is suggested in Evans(2015), which considers the tail probability
\begin{equation}
\Pi_{\Psi}(RB_{\Psi}(\psi\,|\,x)\leq RB_{\Psi}(\psi_{0}\,|\,x)\,|\,x).
\label{strength}%
\end{equation}
(\ref{strength}) can be interpreted as the posterior probability that the true value of $\psi$ has a relative
belief ratio no greater than that of the hypothesized value $\psi_{0}.$ When $RB_{\Psi}(\psi_{0}\,|\,x)<1,$ so there is evidence
against $\psi_{0},$ then a small value for (\ref{strength}) indicates a large
posterior probability that the true value has a relative belief ratio greater
than $RB_{\Psi}(\psi_{0}\,|\,x)$ and so there is strong evidence against
$\psi_{0}.$ When $RB_{\Psi}(\psi_{0}\,|\,x)>1,$ so there is evidence in favor
of $\psi_{0},$ then a large value for (\ref{strength}) indicates a small
posterior probability that the true value has a relative belief ratio greater
than $RB_{\Psi}(\psi_{0}\,|\,x))$ and so there is strong evidence in favor of
$\psi_{0},$ while a small value of (\ref{strength}) only indicates weak
evidence in favor of $\psi_{0}.$

\section{Dirichlet process} \label{Dirichlet}

A relevant summary of the Dirichlet process is presented in this section. The Dirichlet process, formally introduced in Ferguson (1973), is considered the most well-known and widely used prior in Bayesian nonparameteric inference. Specifically, consider  $\mathfrak{X}$ a space with a $\sigma-$algebra $\mathcal{A}$ of subsets of $\mathfrak{X}$. Let $G$ be a fixed probability measure on $(\mathfrak{X},\mathcal{A})$, called the \emph{base measure},  %
 and $a$ be a positive number, called the \emph{concentration parameter}. Following Ferguson (1973), a
random probability measure $P=\left\{  P(A)\right\}  _{A\in\mathcal{A}}$ is
called a Dirichlet process on $(\mathfrak{X},\mathcal{A})$ with parameters $a$
and $G$, denoted by $DP(a,G)$, if for any finite measurable partition $\{A_{1},\ldots,A_{k}\}$ of
$\mathfrak{X}$ with $k \ge 2$, $\left(  P(A_{1}%
),\ldots\,P(A_{k})\right)\sim \text{Dirichlet}(aG(A_{1}),\ldots,$ $aG(A_{k}))$. It is assumed that if
$G(A_{j})=0$, then $P(A_{j})=0$ with a probability one.  For any
$A\in\mathcal{A},$ $P(A) \sim \text{Beta}(aG(A),(1-G(A))$ and so ${E}(P(A))=G(A)\ \ $ and ${Var}(P(A))=G(A)(1-G(A))/(1+a).$
Thus, $G$ plays the role of the center of the process, while $a$ controls concentration, as,  the larger value of $a$, the more likely that $P$ will be close to $G$. Not that, for convenience, we do not distinguish between a
probability measure and its cdf.

An important feature of the Dirichlet process is the conjugacy property. Specifically, if
$x=(x_{1},\ldots,x_{n})$ is a sample from $P\sim DP(a,G)$, then the posterior
distribution of $P$ is $P\,|\,x=P_x\sim DP(a+n,G_{x})$ where
\begin{equation}
G_{x}=a(a+n)^{-1}G+n(a+n)^{-1}F_{n}, \label{DP_posterior}%
\end{equation}
with $F_{n}=n^{-1}\sum_{i=1}^{n}\delta_{{x}_{i}}$ and $\delta_{x_{i}}$ the
Dirac measure at $x_{i}.$ Notice that, $G_{x}$
is a convex combination of the prior base distribution and the empirical
distribution. Clearly,  $G_{x}\to H$  as $a \to \infty$ while   $G_{x}\to F_n$ as $a \to 0$. We refer the reader to Al-Labadi and Zarepour (2013a,b; 2014a) and Al-Labadi and Abdelrazeq (2017) for other interesting asymptotic properties of the Dirichlet process.

Following Ferguson (1973), $P\sim{DP}(a,G)$ has the following series representation
\begin{equation}
P=\sum_{i=1}^{\infty}J_{i} \delta_{Y_{i}}, \label{series-dp}%
\end{equation}
where $\Gamma_{i}=E_{1}+\cdots+E_{i}$, $E_{i} \overset{i.i.d.}\sim \text{exponential}(1)$, $Y_{i} \overset{i.i.d.}\sim  G$
 independent of $\Gamma_{i}$, $L(x)=a\int_{x}^{\infty}t^{-1}e^{-t}dt,x>0,$ $L^{-1}(y)=\inf
\{x>0:L(x)\geq y\}$ and $J_{i}=L^{-1}(\Gamma_{i})/\sum_{i=1}^{\infty
}{L^{-1}(\Gamma_{i})}$.  It follows clearly from (\ref{series-dp}) that a
realization of the Dirichlet process is a discrete probability measure.  This
is correct even when $G$ is absolutely continuous.  Note that, one could resemble the discreteness of $P$ with the discreteness of $F_n$. Since data is always measured to finite accuracy, the true distribution being sampled from is discrete. This makes the  discreteness property of $P$ with no practical significant limitation.  Indeed, by imposing the weak topology, the support for the Dirichlet process
is quite large. Precisely, the support for the Dirichlet process is the set
of all probability measures whose support is contained in the support of the
base measure. This means if the support of the base measure is $\mathfrak{X}$,
then the space of all probability measures is the support of the Dirichlet
process. For instance, if  $G$ is the standard normal, then the Dirichlet
process can choose any probability measure.

Recognizing that no closed form for the inverse of L\'evy measure $L(x)$,  Sethuraman (1994) introduced the stick-breaking approach to define the Dirichlet Process. Specifically, let $(\beta_i)_{i \ge 1}$ be a sequence of i.i.d. random
variables with a $\text{Beta}(1,\alpha)$ distribution. In (\ref{series-dp}), set
\begin{equation}
J_1=\beta_1,~ J_i=\beta_i\prod^{i-1}_{k=1}(1-\beta_k),~i\ge2.\label{eq6}
\end{equation}
and $(Y_{i})_{i \ge {1}}$ independent of $(\beta_i)_{i \ge 1}$. Unlike Ferguson's approach, the stick-breaking construction does not need normalization. By truncating the higher order terms in the sum to  simulate Dirichlet process, we can approximate the Sethuraman stick breaking representation by
\begin{equation}
P_N=\sum_{i=1}^{N}J_{i,N} \delta_{{Y}_i}(\cdot). \label{suth}
\end{equation}
In here,  $(\beta_i)_{i \ge 1}$, $(J_{i,N})_{i \ge 1}$, and $(Y_i)_{i \ge 1}$ are as defined in (\ref{eq6}) with $\beta_N=1$. The assumption that $\beta_N=1$ is necessary to   make the weights add to 1, almost surely (Ishwaran and James, 2001).

For other simulation methods for the Dirichlet process, see
Bondesson (1982),   Wolpert and Ickstadt (1998), Zarepour and Al-Labadi (2012), Al-Labadi and Zarepour (2014b).

\section{Kullback-Leibler Distance}
\label{KL}
Let $F$ and $F_1$ be two continuous cdf's with corresponding probability density functions (pdf's)  $f$ and $f_1$ (with respect to Lebesgue measure). Then Kullback-Leibler divergence or the Relative Entropy between $F$ and $F_1$ is defined as
\begin{eqnarray}
\nonumber d_{KL}(F,F_1)&=&\int_{-\infty}^{\infty} f(x)\log \left(f(x)/f_1(x)\right) dx\\
 &=&-H(F)-\int_{-\infty}^{\infty} f(x) \log f_1(x) dx, \label{KL0}
\end{eqnarray}
where
\begin{equation}
H(F)=-\int_{-\infty}^{\infty} f(x)\log f(x) dx=-E_F\left[\log f(x)\right] \label{entropy}%
\end{equation}
is the \emph{entropy} of $F$ (Shannon, 1948). It is well-know that $d_{KL}(F,F_1) \ge 0$ and the equality holds if and only if $f=f_1$.  However, it is not a distance as it is not symmetric and does not satisfy the triangle inequality (Cover and Thomas, 1991).

The following results provides a simple formula for the distance between $DP(a,G)$  and a continuous cdf. The result uses Al-Labadi, Patel, Vakiloroayaei and  Wan (2018) Bayesian non-parametric estimator of entropy.
\begin{lemma}
\label{KL-distance}
Let $G$ be a continuous cdf with corresponding density function $g$. Let $P_N=\sum_{i=1}^{N}J_{i,N} \delta_{{Y}_i}$ as defined in (\ref{suth}), where  $Y_{1},Y_{2}, $ $\ldots, Y_{N}\overset{i.i.d.}\sim G$ with the probability density function $G^{\prime}(x)=g(x)$. Let $m$ be a positive integer smaller than $N/2$, $Y_{(i)}=Y_{(1)}$ if $i<1$, $Y_{(i)}=Y_{(N)}$ if $i>N$, $Y_{(1)} \le Y_{(2)} \le \cdots \le Y_{(N)}$ are the order statistics of $Y_{1}, Y_{2}, \ldots, Y_{N}$ and
\begin{eqnarray}
H_{m,N,a}=\sum_{i=1}^N J_{i,N} \log \left(\frac{Y_{(i+m)}-Y_{(i-m)}}{c_{i,a}}\right), \label{entropy_est}
\end{eqnarray}
where
$$
c_{i,a} = \left\{
        \begin{array}{ll}
            \sum_{k=2}^{i+m}J_{k,N} & \quad  1\leq i \leq m \\
            \sum_{k=i-m+1}^{i+m}J_{k,N} & \quad  m+1 \leq i \leq N-m\\
            \sum_{k=i-m+1}^{N}J_{k,N} & \quad  N-m+1 \leq i \leq N
        \end{array}
    \right..
$$
Then, as $N \to \infty$, $m \to \infty$, $m/N \to 0$ and $a \to \infty$, we have
\begin{eqnarray*}
d_{KL}(P_{N},G)&=&-H_{m,N,a}-\sum_{i=1}^N J_{i,N} \log \left(g(Y_{(i)})\right)\\
&=&-\sum_{i=1}^N J_{i,N} \log \left\{\frac{\left(Y_{(i+m)}-Y_{(i-m)}\right)g\left(Y_{(i)}\right)}{2mc_{i,a}}\right\} \\
&\overset{p}\to& d_{KL}(P,G),
   \end{eqnarray*}
\end{lemma}

\proof By Lemma 1 of Al-Labadi, Patel, Vakiloroayaei and  Wan (2018), $-H_{m,N, a} \overset{p}\to -H(P)$, where $H(P)$ is the entropy of $P$. Also,
$$\sum_{i=1}^N J_{i,N} \log \left\{g(Y_{(i)})\right\}=\frac{1}{N}\sum_{i=1}^N N J_{i,N} \log  \left\{g(Y_{(i)})\right\}, $$
which, by the weak law of large numbers, converges in probability to
$$E_{P}\left[\log g(Y_{(i)})\right]=\sum_{i=1}^\infty J_{i} \log \left(g(Y_{(i)})\right).$$
Now, by the continuous mapping theorem, we get the result. \endproof

\section{Model Checking Using Relative Belief}

Let $\left\{  F_{\theta}:\theta\in\Theta\right\}  $ denote the collection of continuous
cdf's for the model. Suppose
that $x=(x_{1},\ldots,x_{n})$ is a sample from a distribution $P$. The goal is to assess the hypothesis $\mathcal{H}_{0}:P\in\left\{  F_{\theta}:\theta
\in\Theta\right\}  $. Let the prior on $P$ be $DP(a,G)$ for some choice of $a$ and $G$. Then, by (\ref{DP_posterior}), the posterior distribution is $P\,|\,x\sim DP\left(a+n,G_{x}\right)  $.  As pointed out in the introduction,  if $\mathcal{H}_{0}$ is true, then the posterior distribution of the distance between $P$ and $\left\{ F_{\theta}:\theta\in\Theta\right\}$ should be more concentrated about $0$ than the prior distribution of this distance. So this test will involve a comparison of the concentrations of the prior and posterior distributions of $d_{KL}$ via a relative belief ratio based on $d_{KL}$ with the interpretation as discussed
in$\ $ Section 2. However, to fully implement this approach, it is necessary to discuss the distance measure and the  ideal values for $m,$ $a$ and $G$.

\subsection{Measuring the Distance}

Similar to Al-Labadi and Evans (2018), we compute $d_{KL}(P,F_{\theta(x)})$, where   $F_{\theta(x)}\in$ $\left\{  F_{\theta}:\theta\in\Theta\right\}
$  is the distribution that is best supported by the data. Since the evidence is being measured via relative
belief ratios,  $\theta(x)$ is the relative belief estimate
of $\theta,$ which for the full model parameter is always the same as the  maximum likelihood estimate (MLE).
As such, the value $\theta(x)$ is completely independent of any prior placed
on $\theta.$ Certainly this choice has some asymptotic justification as, under
reasonable conditions, $\theta(x)$ will converge to the best choice (in terms
of Kullback-Leibler divergence) of $\theta$ even when the model fails.

\subsection{The Choice of $m$}
The value of $m$ is required to compute  (\ref{entropy_est}). However, the optimal value $m$  is still an open problem. As discussed in Vasicek (l976), with increasing $N$, the best value of $m$ increases while the ratio $m/N$ tends to zero.  Grzegorzewski and Wieczorkowski (1999) proposed the following formula for optimal values of $m$
 \begin{eqnarray}
\label{optimal-m} m =\floor {\sqrt{N} +0.5},
    \end{eqnarray}
where $\floor {y}$ is the largest integer less than or equal to $y$. Thus, for instance, by (\ref{optimal-m}), if  $N=50$, the best choices of $m$ is 7. In this paper, we will use the rule  (\ref{optimal-m}).
Note that, the value of $m$ in (\ref{optimal-m}) is the value that will be used for the prior. For the posterior, $N$ will be replaced by the number of distinct atoms in $P_N|\,x$,  an approximation of $P|x$. It follows from (\ref{DP_posterior}) that if $a/n$ is close to zero, then the number of distinct  atoms in $P_N|\,x$ will typically be  $n$, the sample size.

\subsection{The Choice of $G$}
Following Al-Labadi and Evans (2018), we set $G = F_{\theta(x)}$ (i.e. $P \sim DP\left(a,F_{\theta(x)}\right)$). There are many benefits of this choice of $G$. First, it avoids prior-data conflict (Evans and Moshonov, 2006; Al-Labadi and Evans, 2017) as the existence of prior-data conflict may lead to the failure having an appreciable concentration of the posterior distribution of $d_{KL} \left( P, F_{\theta(x)} \right) $ about zero, even when $\mathcal{H}_{0}$ is true (Al-Labadi and Evans, 2018). On the other hand, setting $G =  F_{\theta(x)} $ would appear to induce a data dependent prior distribution for $d_{KL}$. The following lemma implies that this is not the case and so, with this choice, the approach is prior distribution-free.
\begin{lemma} \label{ind} If $P \sim DP\left(a,F_{\theta(x)}\right)$, then  the distribution of $d_{KL}\left(P,F_{\theta(x)}\right)$ does not depend on $F_{\theta(x)}$.
\smallskip
\end{lemma}
\proof By Lemma \ref{KL-distance},
 \begin{eqnarray*}
d_{KL}(P_{N},F_{\theta(x)})&=&-\sum_{i=1}^N J_{i,N} \log \bigg(\frac{1}{{c_{i,a}}}\frac{Y_{(i+m)}-Y_{(i-m)}}
{F_{\theta(x)}\left(Y_{(i+m)}\right)-F_{\theta(x)}\left(Y_{(i-m)}\right)} \times \\
&&{\left[F_{\theta(x)}\left(Y_{(i+m)}\right)-F_{\theta(x)}\left(Y_{(i-m)}\right)\right]} f_{\theta(x)}(Y_{(i)})\bigg)
   \end{eqnarray*}
Note that, as $m \to \infty$ such that $m/N \to 0$, we have
 \begin{eqnarray}
\label{derivative} \frac{Y_{(i+m)}-Y_{(i-m)}}{F_{\theta(x)}\left(Y_{(i+m)}\right)-F_{\theta(x)}\left(Y_{(i-m)}\right)}&=&\frac{1}{f_{\theta(x)}(Y_{(i)})},
    \end{eqnarray}
where $f_{\theta(x)}$ is the pdf of $F_{\theta(x)}$.
Also, since $(Y_i)_{i \ge 1}$ is a sequence of
i.i.d. random variables with continuous distribution $F_{\theta(x)},$ for $i \ge 1,$ we have $U_i\overset{d}=F_{\theta(x)}(Y_i)$, where $\left(U_i\right)_{i \ge 1}$ is a sequence of i.i.d. random variables with a uniform distribution on $[0,1]$. Thus,
 \begin{eqnarray}
d_{KL}(P_{N},F_{\theta(x)})&\overset{d}=&-\sum_{i=1}^N J_{i,N} \log \bigg(\frac{U_{(i+m)}-U_{(i-m)}}{c_{i,a}}\bigg). \label{uniform}
   \end{eqnarray}
Now, as $N \to \infty$, $m \to \infty$, $m/N \to 0$, by Lemma \ref{KL-distance}, we conclude that the distribution of $d_{KL}(P,F_{\theta(x)})$ does not depend on $F_{\theta(x)}$.
\endproof

Note that, similar to Noughabi and Arghami (2013), if $G(y)=y+b$, $G(y)=ay$ or $G(y)=ay+b$, which involve the case of location, scale and location-scale families, then
(\ref{derivative}) holds without any condition.

The following result shows that the posterior distribution of $d_{KL}\left( P, F_{\theta\left(x\right)}\right)$ becomes concentrated around 0 as sample size increases if and only if $\mathcal{H}_{0}$ is true. The proof follows straightforwardly from the properties of the KL divergence and (\ref{DP_posterior}). Thus it is omitted.

\begin{lemma} Let $x=(x_1,\ldots,x_n)\sim P$, where $P \sim DP\left(a,F_{\theta(x)}\right)$. Suppose that
$\theta(x)\overset{a.s.}{\rightarrow}\theta_{0},\sup_{y}|F_{\theta
(x)}(y)-F_{\theta_{0}}(y)|\overset{a.s.}{\rightarrow}0$ as $n\rightarrow
\infty.$
\begin{enumerate}
\item [(i)] If $\mathcal{H}_{0}$ is true, then, as $n\to \infty$, $d_{KL}\left(  P|x,F_{\theta
(x)}\right)  \overset{a.s.}{\rightarrow}0$.

\item [(ii)] If $\mathcal{H}_{0}$\ is
false, then, as $n\to \infty$, $\lim\inf d_{KL}(P|x,F_{\theta(x)})\overset{a.s.}{>}0.$
\end{enumerate}
\end{lemma}

\subsection{The Choice of $a$}

The selection of $a$ is very important. In principle, larger values of $a$ must be chosen to detect smaller deviations. Therefore, it is possible to consider several values of $a$. For example, one may start with $a = 1$. If, as the
value of $a$ is increased, the corresponding relative belief ratio drops
rapidly below 1, then this is a clear indication against $\mathcal{H}_{0}$. As
will be seen in the examples, when the null hypothesis is correct, the relative belief
ratio always remains above 1 when larger values of $a$ are considered. On the other hand, if the relative belief ratio is below
than 1 and, as the value of $a$ is increased (i.e., using a more concentrated prior), the corresponding relative belief ratio increases above 1, then this is a good  indication in favour of $\mathcal{H}_{0}$. It is highly recommended to choose $a\leq0.5n$,  however, otherwise the prior may become too influential. See Al-Labadi and Zarepour (2017) for the justification of this recommendation. It is noticed that, for most purposes, setting $a$ between $1$ and $10$ is found satisfactory. This choice of $a$  is also recommended by Holmes, Caron, Griffin and Stephens (2015) when using the P\'{o}lya tree prior for the two-sample problem.  This issue is further discussed in  Table \ref{tab1} of Section 7.

The following result is useful in the elicitation process of $a$.

\begin{lemma}  If $P\sim DP(a,F_{\theta(x)}),$ then
\begin{eqnarray}
\nonumber		E\left[d_{KL}(P_{N},F_{\theta(x)})\right] &=& \frac{2}{N}\sum_{i=1}^{m}\left[\psi\left(\frac{a(m+i-1)}{N}+1\right)-\psi(m+i-1)\right] \\
\nonumber&&+\frac{N-2m}{N} \left(\psi\left(\frac{2am}{N}+1\right)-\psi(2m)\right)\\
&&+\psi(N+1)-\psi(a+1) ,\label{expect}
\end{eqnarray}
 where $\psi(x) = \Gamma'(x)/\ \Gamma(x)$ is the digamma function.
\end{lemma}
\proof
By  (\ref{uniform}) and independence,
 \begin{eqnarray}
\nonumber E\left[d_{KL}(P_{N},F_{\theta(x)})\right]&=& -\sum_{i=1}^N E\left[J_{i,N}\right] E\left[\log \left(U_{(i+m)}-U_{(i-m)}\right)\right]\\
&&+ \sum_{i=1}^N E\left[J_{i,N} \log c_{i,a}\right]. \label{expectation0}
   \end{eqnarray}
Since $U(i)=U(1)$ for $i <1 $ and $U(i)=U(N)$ for $i>N$ and using the well-known fact that $U_{(s)}-U_{(r)} \sim Beta(s-r,N-s+r+1)$ we have:
\begin{eqnarray}
	\nonumber	\lefteqn{\sum_{i=1}^{N} E\left[\log(y_{(i+m)} - y_{(i-m)})\right]}\\
\nonumber &=&\sum_{i=1}^{m} E\left[\log(U_{(i+m)} - U_{(i-m)})\right]+ \\
\nonumber &&\sum_{i=m+1}^{N-m} E\left[\log(U_{(i+m)} - U_{(1)})\right]+\sum_{N-m+1}^N E\left[\log(U_{(N)} - U_{(i-m)})\right]\\
	\nonumber	 &=& (N-2m) \left(\psi(2m) - \psi(N+1)\right)+ 2 \sum_{i=1}^{m} \left(\psi(i+m-1)-\psi(N+1)\right)\\
		  &=& (N-2m) \psi(2m)-N \psi(N+1) + 2 \sum_{i=1}^{m} \psi(i+m-1). \label{expectation1}
\end{eqnarray}
On the other hand,
\begin{eqnarray}
	\nonumber	\lefteqn{\sum_{i=1}^{N} E\left[J_{i,N} \log c_{i,a}\right]}\\ \label{expectationI} &=&\sum_{i=1}^{m} E\left[J_{i,N}\log \left(\sum_{k=2}^{i+m}J_{k,N}\right)\right]+ \\
\label{expectationII} &&\sum_{i=m+1}^{N-m}E\left[J_{i,N}\log \left(\sum_{k=i-m+1}^{i+m}J_{k,N}\right)\right]+\\
  \label{expectationIII}&&\sum_{i=N-m+1}^{N} E\left[J_{i,N}\log \left(\sum_{k=i-m+1}^{N}J_{k,N}\right)\right].
\end{eqnarray}
From the proof of Lemma 1  of Al-Labadi, Patel, Vakiloroayaei and  Wan (2018) and using the facts that $\psi(x+1)=\psi(x)+1/x$  (Abramowitz and Stegun, 1972), we have
\begin{eqnarray*}
		(\ref{expectationI})&=&\frac{1}{N}\sum_{i=1}^{m} \psi\left(\frac{a(m+i-1)}{N}+1\right)-\frac{m}{N}\psi(a+1),
\end{eqnarray*}
\begin{eqnarray*}
		(\ref{expectationII})&=&\frac{1}{N}\sum_{i=m+1}^{N-m}\psi\left(\frac{2am}{N}+1\right)-\frac{N-2m}{N}\psi(a+1)\\
                             &=&\frac{N-2m}{N}\psi\left(\frac{2am}{N}+1\right)-\frac{N-2m}{N}\psi(a+1)
\end{eqnarray*}
and
\begin{eqnarray*}
		(\ref{expectationIII})&=&\frac{1}{N}\sum_{i=N-m+1}^{N} \psi\left(\frac{a(N+m-i)}{N}+1\right)-\frac{m}{N}\psi(a+1)\\
                              &=&\frac{1}{N}\sum_{i=1}^{m} \psi\left(\frac{a(m+i-1)}{N}+1\right)-\frac{m}{N}\psi(a+1).
\end{eqnarray*}
Substitute (\ref{expectation1}), (\ref{expectationI}), (\ref{expectationII}) and (\ref{expectationIII}) in (\ref{expectation0}), we get the result.
\endproof

%
%
%

\section{Computational Algorithm}

To use (\ref{relbel}), closed forms of the prior and posterior densities of $D=d_{KL}(P,F_{\theta
(x)})$ are required, which is typically not available. Consequently, the relative belief ratio needs to
be approximated via simulation. A particular attention here should  be given to the case when both $\pi_{D}(0\,|\,x) $ and $\pi _{D}(0) $ are close to $0$.  In such a case, determining $RB_{D}(0\,|\,x)$  is challenging. However, as discussed in Section 2, the formal definition of $RB_{D}(0\,|\,x)$ is given as a limit and this limit can be approximated by $RB_{D}([0,d_{\ast})\,|\,x)$, the ratio of the posterior to prior probability that $0\leq D\leq d_{\ast},$ for a suitably small value of $d_{\ast}.$

We adapt the procedure outlined in Al-Labadi and Evans (2018). This approach is based on $M$
quantiles of the prior distribution of $D,$ namely, the $i$-th interval is
$[d_{i/M}(pr),d_{(i+1)/M}(pr))$ where $d_{i/M}(pr)$ is the $(i/M)$-th quantile for
$i=0,\ldots,M.$ Note that values in the left tail of this distribution
correspond to those $P,$ in the population of distributions, that, according
to the prior at least, do not differ materially from 0. As such we will consider the left-tail quantile of this prior distribution, such as the
$0.05$-quantile or the $0.01$-quantile, so $d_{\ast}=d_{i_{0}/M}(pr)$ where
$i_{0}/M\approx0.05$ or $i_{0}/M\approx0.01.$

The following gives a computational algorithm for assessing $\mathcal{H}_{0}$.

\noindent\textbf{Algorithm A }\textit{(}\emph{Relative belief algorithm for
model checking}\textit{):\\}
\noindent (i)  Generate a sample from $P_N$, where $P_N $ is  an approximation of    $P \sim DP(a,F_{\theta
(x)})$. See Section 3. \\
\noindent (ii) Compute $d(pr)=d_{KL}(P_N,F_{\theta(x)})$.\\
\noindent \noindent (iii) Repeat steps (i) and (ii) to obtain a sample of $r_{1}$ values from the
prior of $D$.  \\
\noindent (iv) Generate  a sample from $P_N|x$, where $P_N|x$ an approximation of   $P|x \sim DP(a+n,G_{x}%
)$. \\
\noindent (v) Compute $d(po)=d_{KL}(P_N|x,F_{\theta(x)})$. \\
\noindent (vi) Repeat steps (iv)-(v) to obtain a sample of $r_{2}$ values from the
posterior of $D$.\\
\noindent (vii) For a fixed positive number $M$, let $\hat{F}_{D}$ denote the empirical cdf
of $D$ based on the prior sample in (iii) and for $i=0,\ldots,M,$ let $\hat
{d}_{i/M}(pr)$ be the estimate of $d_{i/M}(pr),$ the $(i/M)$-th prior quantile of $D.$
Here $\hat{d}_{0}(pr)=0$, and $\hat{d}_{1}(pr)$ is the largest value of $d(pr)$. Let
$\hat{F}_{D}(\cdot\,|\,x)$ denote the empirical cdf of $D$ based on the
posterior sample of $d(po)$ in (vi). For $d\in\lbrack\hat{d}_{i/M}(pr),\hat{d}_{(i+1)/M}(pr))$,
estimate $RB_{D}(d\,|\,x)$ by the ratio of the estimates of the posterior and prior contents of $[\hat
{d}_{i/M}(pr),\hat{d}_{(i+1)/M}(pr)).$ Specifically,
\begin{equation}
\widehat{RB}_{D}(d\,|\,x)=M\{\hat{F}_{D}(\hat{d}_{(i+1)/M}(pr)\,|\,x)-\hat{F}%
_{D}(\hat{d}_{i/M}(pr)\,|\,x)\},\label{rbest}%
\end{equation}
 Moreover, estimate $RB_{D}(0\,|\,x)$ by
$\widehat{RB}_{D}(0\,|\,x)=M\widehat{F}_{D}(\hat{d}_{i_{0}/M}(pr)\,|\,x)$ where
$i_{0}$ is chosen so that $i_{0}/M$ is not too small (typically $i_{0}%
/M\approx0.05)$.\\
\noindent (viii) Estimate the strength $DP_{D}(RB_{D}(d\,|\,x)\leq RB_{D}%
(0\,|\,x)\,|\,x)$ by the finite sum
\begin{equation}
\sum_{\{i\geq i_{0}:\widehat{RB}_{D}(\hat{d}_{i/M}(pr)\,|\,x)\leq\widehat{RB}%
_{D}(0\,|\,x)\}}(\hat{F}_{D}(\hat{d}_{(i+1)/M}(pr)\,|\,x)-\hat{F}_{D}(\hat
{d}_{i/M}(pr)\,|\,x)). \label{strest}%
\end{equation}

For fixed $M,$ as $r_{1}\rightarrow\infty,r_{2}\rightarrow\infty,$ then
$\hat{d}_{i/M}(pr)$ converges almost surely to $d_{i/M}(pr)$, (\ref{rbest}) converge almost surely to $RB_{D}(d\,|\,x)$ and
(\ref{strest}) converge almost surely to $DP_{D}%
(RB_{D}(d\,|\,x)\leq RB_{D}(0\,|\,x)\,|\,x)$ (Al-Labadi and Evans, 2018).

\section{Examples}

In this section, the approach is illustrated through three examples. In all the
examples, the prior was taken to be $DP\left(  a,F_{\theta(x)}\right)$ and, in Algorithm A,  we set $r_{1}=r_{2}=2000$, $N=200$, $M=20$ and $i_{0}=1$. A critical factor here for success are the choices of $a$   as the
prior has to be sufficiently concentrated about the family. The sensitivity to the choice of $\ a$  is  investigated and we record only a few values in the tables.

\noindent\textbf{Example 1. }\textit{Location normal model.}

In this example, samples of $n=20$ was
generated from the distribution $P_{true}$ in  Table \ref{tab1}. Then the methodology was
applied to assess whether or not the correct model is $\{F_{\theta}:\theta\in\Theta\}=\{N(\theta,1):\theta
\in\mathcal{%
\mathbb{R}
\}}$ and so $\theta(x)=\bar{x}.$ Thus, by Lemma \ref{KL-distance},
\begin{eqnarray*}
d_{KL}(P_{N},F_{\theta(x)})&=&-H_{m,N,a}-\sum_{i=1}^N J_{i,N} \log \left(f_{\theta(x)}(Y_{(i)})\right),
   \end{eqnarray*}
where $$f_{\theta(x)}(Y_{(i)})=\frac{1}{\sqrt{2\pi}}e^{-\frac{1}{2}\left(Y_{(i)}-\bar{x}\right)^2}.$$
It follow that
\begin{eqnarray*}
d_{KL}(P_{N},F_{\theta(x)})&=&-H_{m,N,a}+\frac{1}{2}\log (2\pi)+\frac{1}{2}\sum_{i=1}^N J_{i,N} \left(Y_{(i)}-\bar{x}\right)^2.
   \end{eqnarray*}
In Table \ref{tab1} the relative belief
ratios and the strengths are recorded for testing the location normal model
against a variety of alternatives with two choices of the hyperparameter $a$ and $m$. Recalling that we want $RB>1$ and the strength close to 1 when $\mathcal{H}%
_{0}$ is true and $RB<1$ and the strength close to 0 when $\mathcal{H}_{0}$ is
false, it is seen that the methodology using $d_{KL}(P,F_{\theta(x)})$
performs well in every instance.

\begin{table}[h] \centering
\begin{tabular}
[c]{ccccc}\hline
$P_{true}$ & $a$ & $d_{0.05}(pr)$ &$RB$ & Strength  \\\hline

$N(0,1)$ &1&0.5689 & 20    &    1\\
&  $5$ &0.1441 &13.2041& 0.3395 \\
&  $10$ &0.0546& 4.4124& 0.7780  \\

$N(10,1)$ &$1$  &0.5690 & 20 &       1  \\
&  $5$ &0.1440684 &13.2041 &0.3394  \\
&  $10$ &0.0546& 4.4124 &0.7779 \\

$N(0,4)$  &$1$ &0.5512 &0.8221&   0.1355\\
&  $5$ & 0.1283 &0.03988  & 0.0015\\
&  $10$ &0.0494& 0.0300&    0.0000\\

$0.5N(-2,1)+0.5N(2,1)$ &$1$& 0.5440& 0.3878  &      0  \\
&  $5$ &0.1286 &  0  &      0 \\
&  $10$ &0.0548 & 0  &      0 \\

$t_{0.5}$ &$1$&0.570 & 0  &      0 \\
&  $5$ &0.1331  & 0  &      0  \\
&  $10$ &0.0569  & 0 &       0  \\


$t_3$ &$1$&0.5485 &4.6007& 0.9980 \\
&  $5$ &0.1284 &0.8618 &0.3577  \\
&  $10$ &0.0511 &0.8819 &0.5147\\

Cauchy$(0,1)$ &$1$&0.5390 &0.0099&        0  \\
&  $5$ &0.1441 &  0  &      0 \\
&  $10$ &0.0548  & 0 &       0 \\
&  &  &
\end{tabular}
\caption{Relative belief ratios and strengths for testing the location normal model with various alternatives and choices of
$a$  in Example 1.}\label{tab1}%
\end{table}%

%
%

\newpage

\noindent\textbf{Example 2.} \textit{The Gumbel Model.}

In this example, we consider the Gumbel model. This model is commonly used  in environmental sciences, hydrology in the modeling of heavy rain, floods and industrial applications. A random
variable $Y$ is said to have the Gumbel distribution if its
probability density function has the form
$$f(y;\xi,\beta)=\frac{1}{\beta}\exp\left\{ -\frac{y-\xi}{\beta}-\exp\left(-\frac{y-\xi}{\beta}\right)\right\}, \ \  y, \xi  \in \mathbb{R},\ \ \beta>0.$$
Here $\xi$ represents the location parameter and $\beta$
represents the scale parameter. The
following dataset gives the annual maxima of daily rainfall (in mm) during the period 1967-2001 recorded at the \'Alamo, Veracruz, meteorological station, M\'exico maximum flood levels of the Susquehanna River at
Harrisburg, Pennsylvania, over four-year periods (1890-1969) in millions of
cubic feet per second.
\begin{quote}
$86.8, 78.5, 93.1, 95.5, 78.1, 89.9, 109.5, 161.6, 187.6, 89.9, 73.4, 78.1,  73.3, \newline
  130.1, 188.3, 113.9, 42.5, 80.0, 142.6, 42.9, 60.2, 100.0, 129.0,  98.0,  116.4, \newline
 37.9, 60.7, 48.7, 39.7, 80.3, 30.7, 120.0, 160.0, 64.3, 80.0$

\end{quote}
\noindent According to  P\'erez-Rodr\'iguez, Vaquera-Huerta and Villase\~nor-Alva (2009), the maximum
likelihood estimators for $\xi$ and $\beta$ are $74.5432 $ and $32.4328$, respectively. The goal is to test whether the underlying
distribution is Gumbel. The results in Table \ref{tab3} indicate that indeed the data can be considered as coming from a
Gumbel distribution as there is evidence in favor of this model.%

\begin{table}[h] \centering
\begin{tabular}
[c]{cccccc}\hline
$a$ & 1 & 5 & 10 & 15 & 20\\
$d_{0.05}(pr)$ & $0.5573$ & $0.1209$ & $0.0499$ & $0.0306$ &  $0.0215$ \\
$RB$ & $20$ & $13.2132$ & $5.7211$ & $3.7904$ & $3.0154$ \\
Strength & $1$ & $1$ & $1$ & $1$ & $1$ \\
\hline
& & & & &
\end{tabular}
\caption{Relative belief ratios and strengths for testing the Gumbel model in Example 2.}\label{tab3}%
\end{table}%

\bigskip

\noindent\textbf{Example 3.} \textit{Lifetimes of Kevlar pressure vessels.}

Consider  the data of 100 stress-rupture lifetimes of Kevlar pressure
vessels presented in Andrews and Herzberg (1985). The goal is to test whether the underlying distribution is normal. That is, $\{F_{\theta}:\theta\in\Theta\}=\{N(\mu,\sigma^{2}%
):\theta=(\mu,\sigma^{2})\in\mathcal{%
\mathbb{R}
}\times(0,\infty)\}$ and so $\theta(x)=(\bar{x},\sum_{i=1}^{n}(x-\bar{x}%
)^{2}/n).$ For this data set, $\theta(x)=(209.171,37606.56)$. Previous studies such as Evans and Swartz
(1994), Verdinelli and Wasserman (1998) and Al-Labadi and Evans (2018) suggested that model is not
correct.  The results in Table \ref{tab4} support the non-normality of this data set only
when using a more concentrated prior.

\begin{table}[h] \centering
\begin{tabular}
[c]{cccccccc}\hline
$a$ & 1 & 5 & 10 & 15 & 20 & 25&30\\
$d_{0.05}(pr)$ &0.5451 & 0.1194& 0.0501& 0.0316&  0.0240 &0.0146 & 0.0130 \\
$RB$ & 17.65 &  1.8107 & 0.6847  & 0.4782&0.3915& 0.1784&0.1697 \\
Strength & 1 &0.7009 &   0.0342 &  0.7009 &    0.1037 &0.0355 &0.0089\\
\hline
& & & & & & &
\end{tabular}
\caption{Relative belief ratios and strengths for testing the normality of the Kevlar data and various choices of
$a$ in Example 3.}\label{tab4}%
\end{table}%

\section{Conclusions}

A general procedure for model checking based on integrating the Dirichlet process, the Kullback-Leibler divergence and the relative belief ratio has been considered.  The offered approach solved the issue that Dirichlet process is a discrete probability measure with probability 1 and the Kullback-Leibler divergence between any discrete distribution and any continuous distribution is infinity. Applying the approach is simple and does not require obtaining a closed form of the relative belief ratio.  Numerous examples are presented in which the proposed approach shows excellent performance.

\end{document}